\documentclass[12pt]{article}
\usepackage{amssymb,latexsym,amsmath}
\usepackage{graphics}                 
\usepackage{color}                    
\usepackage{hyperref}                 
\usepackage[all]{xy}

\begin{document}
\pagestyle{plain} \headheight=5mm \topmargin=-5mm 

\newtheorem{definition}{Definition}
\newtheorem{theorem}{Theorem}[section]
\newtheorem{proposition}[theorem]{Proposition}
\newtheorem{lemma}[theorem]{Lemma}
\newtheorem{corollary}[theorem]{Corollary}
\newtheorem{question}[theorem]{Question}
\newtheorem{remark}[theorem]{Remark}
\newtheorem{example}[theorem]{Example}

\newtheorem{correction}{Correction}

\def\bp{{\bf Proof.}\hspace{2mm}}
\def\qe{\hfill$\Box$}

\def\A{{\mathbb{A}}}
\def\C{{\mathbb{C}}}

\def\L{{\mathbb{L}}}
\def\P{{\mathbb{P}}}
\def\Q{{\mathbb{Q}}}
\def\Z{{\mathbb{Z}}}
\def\Ch{{\rm Ch}}
\def\p{{\mathbf{p}}}

\def\sp{{\rm SP}}
\def\rank{{\rm rank}}

\def\mC{{\mathcal{C}}}
\def\D{{\mathcal{D}}}

\title{Algebraic Cycles of a Fixed Degree}
\author{Wenchuan Hu}
\maketitle

\begin{abstract}
In this paper, the homotopy groups of  Chow variety $C_{p,d}(\P^n)$ of effective $p$-cycles of degree $d$ is proved to be stable in the sense
that  $p$ or $n$  increases. We also obtain a  negative answer  to a question  by Lawson and Michelsohn
on  homotopy groups for the space of degree two cycles.
\end{abstract}


\section{Introduction}
Let $\P^n$ be the  complex projective space of  dimension $n$ and let $C_{p,d}(\P^n)$
be the space of effective algebraic $p$-cycles of degree $d$ on $\P^n$. A
fact proved by Chow and Van der Waerden is that $C_{p,d}(\P^n)$ carries the structure of a closed complex algebraic variety \cite{Chow-Waerden}.
Hence it carries the structure of a compact Hausdorff space.

Consider the spaces
$$\D(d):=\lim_{p,q\to \infty} C_{p,d}(\P^{p+q})$$
of a fixed degree (with arbitrary dimension and codimension), as introduced in \cite{LM}.  Then there is a filtration
$$BU=\D(1)\subset\D(2)\subset\cdots\subset \D(\infty)=K(even,\Z),
$$
where $BU=\lim_{q\to \infty}BU_q$ and  $K(even,\Z)=\prod_{i=1}^{\infty} K(2i,\Z)$ (the weak product of
Eilenberg -Maclane spaces).

The inclusion map $\D(d)\subset \D(\infty)$ induces maps on homology and homotopy groups.
It was proved in \cite{LM} that $\D(1)\subset \D(\infty)$ induces an injective on homotopy groups.
 Moreover, as abstract groups $\pi_*(\D(1))\cong \pi_*(\D(\infty))$.

The following natural question was proposed in \cite{LM}:
\begin{question}\label{ques1}
 Is $\pi_*(\D(d))\to \pi_*(\D(\infty))$ injective for $d\geq 1$?
\end{question}

An affirmative answer to Question \ref{ques1} implies that $\pi_*(\D(d))\cong \pi_*(\D(\infty))$ as abstract groups.

The  first main result in this paper is the following  negative answer to Question \ref{ques1} for $d=2$.
\begin{theorem}\label{Th1.2}
There is an integer $k>0$ such that the induced map $\pi_k(\D(2))\to \pi_k(\D(\infty))$ from the inclusion
$\D(2)\subset\D(\infty)$ is not injective. Moreover, there is an  integer
 $k>0$ such that, as abstract groups, $\pi_k(\D(2))\ncong \pi_k(\D(\infty))$.
\end{theorem}

The proof of this theorem is based on   Theorem \ref{prop1.1} below and direct calculations under the assumption of a positive answer to Question \ref{ques1}.

The second main result is the following:

\begin{theorem}\label{prop1.1}
 $\pi_k(\D(d))\cong \pi_k(C_{p,d}(\P^n))$ for $k\leq min\{2p+1, 2(n-p)\}$.
\end{theorem}

The proof of Theorem \ref{prop1.1} is given in section \ref{sec:3}.

\medskip
\thanks{\emph{Acknowledgements}: I would like to gratitude
Eric Friedlander and Blaine Lawson for their interest 
and helpful advice on the organization of the paper.
}
\section{Homology groups of the space of algebraic cycles with degree two}

The $p$-cycles of degree 1 is linear, so $C_{p,1}(\P^n)=G(p+1,n+1)$,
where $G(p+1,n+1)$ denotes the Grassmannian of $(p+1)$-plane in $\C^{n+1}$. Note that
\begin{equation}\label{eq1}
C_{p,2}(\P^n)=\sp^2(G(p+1,n+1))\cup Q_{p,n}
\end{equation}
where $\sp^i(X)$ denotes the $i$-th symmetric product $X$ and $Q_{p,n}$ consists of
effective irreducible $p$-cycles of degree 2 in $\P^n$. An irreducible $p$-cycle of degree 2 always degenerates if $n\geq p+2$
(cf. \cite{Griffiths-Harris}).
 Hence $$Q_{p,n}=\big\{(\P^{p+1},c)\in G(p+2,n+1)\times C_{p,2}(\P^{n})|\hbox{$c$ is irreducible and
 $c\in C_{p,2}(\P^{p+1})$} \big\}$$
 and it is a fiber bundle over $G(p+2,n+1)$ whose fiber is the space of all irreducible quadric hypersurfaces
 in $\P^{p+1}$.
 Note that all irreducible quadric hypersurfaces in $\P^{p+1}$ is isomorphic
 to $\P^{(^{p+3}_{~2})-1}-\sp^2(\P^{p+1})$, i.e., the complement of non-irreducible
  quadrics (which is a pair of $p$-planes) in the space of all quadric hypersurfaces in $\P^{p+1}$.

To prove Theorem \ref{Th1.2}, we assume that the answer to  Question \ref{ques1} is
affirmative for $d=2$. Then $\pi_{2k}(\D(2))$ is a subgroup of $\Z$ and so $\pi_{2k}(\D(2))\cong 0$ or $\Z$.
 Note that the map
$\pi_{2k}(\D(1))\to \pi_{2k}(\D(\infty))\cong \Z$ is multiplication by $(k-1)!$ (cf. \cite{LM},
Theorem 4.4) and it factors through $\pi_{2k}(\D(2))$. So  $\pi_{2k}(\D(2))$ is nontrivial and $\pi_{2k}(\D(2))\cong \Z$
for all $k$ if Question \ref{ques1} have a positive answer. Similarly, $\pi_{2k-1}(\D(2))=0$ by assuming the positive answer of Question \ref{ques1}.

By Theorem \ref{prop1.1}, we have
$$\pi_k(C_{p,2}(\P^n))\cong\pi_k(\D(2))=\left\{
\begin{array}{ll}
\Z, & \hbox{$k\leq 2p+1$ and even;}\\
0,  & \hbox{$k\leq 2p+1$ and odd.}
\end{array}
\right.
$$

In the following computation, we take $p=4,d=2$ as our example.

\begin{lemma}\label{lemma2.1}
Let $X \to B$ be a fibration between CW complexes with fiber $F$.
Suppose that $B$ is simply connected, $H_{2i}(F,\Q)$ is  finite dimensional, $H_{2i-1}(B,\Q)$ and $H_{2i-1}(F,\Q)$ vanish.
Then $H_k(X,\Q)\cong \bigoplus_{i+j=k} H_i(B,\Q)\otimes H_j(F,\Q)$, that is, the Leray spectral sequence degenerates at $E^2$.
\end{lemma}
\bp
By the  Leray's Theorem for singular homology, we get the $E^2$ term
$$
E^2_{p,q}=H_p(B, H_q(F,\Q))\cong H_p(B,\Q)\otimes H_q(F,\Q), \quad d^2:E^2_{p,q}\to E^2_{p-2,q+1}
$$
since $B$ is simply connected.

From the assumption,  all odd dimensional homology groups of  $B$
and $F$ vanish, so at least one of $E^2_{p,q}$ and $E^2_{p-3,q+2}$ vanishes.
This implies that $d^2$ is a zero map . Hence we get $E^3_{p,q}=E^2_{p,q}$ and
$d^3:E^2_{p,q}\to E^2_{p-3,q+2}$. By the same reason, $d^3=d^4=\cdots=0$. Therefore, the Leray spectral sequence
degenerates at $E^2$, i.e.,
$$\bigoplus_{p+q=k}H_p(B,\Q)\otimes H_q(F,\Q)\cong\bigoplus_{p+q=k}E^2_{p,q}=\bigoplus_{p+q=k} E^{\infty}=H_k(X,\Q).$$

\qe

\begin{proposition}\label{prop2.2}
 Let $X$ be a connected CW complex such that
 $$\pi_k(X)=\left\{
 \begin{array}{ll}
 \Z, & \hbox{$0<k\leq 9$ and even;}\\
0,  & \hbox{$k\leq 9$ and odd.}
 \end{array}
  \right. $$
 Then the first 10 Betti numbers  $\beta_i(X)$ of $X$  are
\begin{equation}
\beta_i(X)=\left\{
\begin{array}{lll}
 1,1,2,3,5,&\hbox{for $i=0,2,4,6,8$}\\
 0, & \hbox{for $i=1,3,5,7,9$}
\end{array}
 \right.
\end{equation}
\end{proposition}
\bp Let $\cdots\to Y_n\to Y_{n-1}\to\cdots \to Y_1=K(\pi_1(X),1)$ be the Postnikov approximation
of $X$, where $Y_n\to Y_{n-1}$ is a fibration with $K(\pi_n(X),n)$ as fibers (cf. \cite{Whitehead}, Chapter IX).
For a fixed $n$, we have isomorphisms of  homotopy groups $\pi_q(X)\cong \pi_q(Y_n)$ and homology groups
$H_q(X,\Q)\cong H_q(Y_n,\Q)$
for $1\leq q\leq n$. Therefore, the first 10 Betti numbers of $X$ coincide with those of $Y_9$.

Note that $Y_2$ is
homotopy equivalent (denote by $\simeq$) to $K(\Z,2)$ since $Y_1$ is contractible. Since $Y_3\to Y_2$ is
a fibration with $K(\pi_3(X),3)\simeq *$ as fibers, we get $Y_3\simeq Y_2$.
Note that $Y_4\to Y_3$ is a fibration with $K(\pi_4(X),4)=K(\Z,4)$ as fibers,
we obtain $H_k(X,\Q)\cong \bigoplus_{i+j=k} H_i(Y_3,\Q)\otimes H_j(K(\Z,4),\Q)$
by Lemma \ref{lemma2.1}. Using Lemma \ref{lemma2.1} for several times,
we get (modulo $H_*(-,\Q)$ for $*>9$)
$$\begin{array}{lll}
 H_*(X,\Q)&\cong& H_*(Y_9,\Q)\\
 &\cong& H_*(Y_8,\Q), \hbox{(since $Y_9\simeq Y_8$)}\\
  &\cong& H_*(Y_7,\Q)\otimes H_*(K(\Z,8),\Q),\hbox{(since $K(\Z,8)\to Y_8\to Y_7$ is a fibration)}\\
 &\cong& H_*(Y_6,\Q)\otimes H_*(K(\Z,8),\Q)\\
 &\cong & H_*(Y_4,\Q)\otimes H_*(K(\Z,6),\Q)\otimes H_*(K(\Z,8),\Q)\\
 &\cong & H_*(Y_2,\Q)\otimes  H_*(K(\Z,4),\Q)\otimes H_*(K(\Z,6),\Q)\otimes H_*(K(\Z,8),\Q)\\
 &\cong & H_*(K(\Z,2),\Q)\otimes  H_*(K(\Z,4),\Q)\otimes H_*(K(\Z,6),\Q)\otimes H_*(K(\Z,8),\Q)
\end{array}
$$

Therefore, the first 10 Betti numbers $\beta_i(X)$ of $X$ are given as follows:
$$
\beta_i(X)=\left\{
\begin{array}{lll}
 1,1,2,3,5;&\hbox{for $i=0,2,4,6,8$}\\
 0; & \hbox{for $i=1,3,5,7,9$.}
\end{array}
 \right.
$$

\qe

The combination of Theorem \ref{prop1.1} and Proposition \ref{prop2.2} implies the following result:

\begin{corollary}\label{cor2.1}
 If the answer to Question \ref{ques1} is affirmative for $d=2$, then the first 10 Betti numbers of
 $C_{4,2} (\P^n)$ ($n\geq 9$) are given by
 $$
\beta_i(C_{4,2} (\P^n))=\left\{
\begin{array}{lll}
 1,1,2,3,5;&\hbox{for $i=0,2,4,6,8$}\\
 0; & \hbox{for $i=1,3,5,7,9$.}
\end{array}
 \right.
$$

\end{corollary}

From the proof of Proposition \ref{prop2.2}, we actually prove the following result:
\begin{remark}
Let $M$ be a connected CW complex such that $\pi_k(X)=0$ for $k$
odd and $\pi_k(M)\cong \Z$ for $k$ positive even integers. Then
$$
\rank (H_k(M))=\left\{
\begin{array}{lll}
p(k), &\hbox{if $k$ is even}\\
0,& \hbox{if $k$ is odd}
\end{array}
\right.
$$
where $p(k)$ represents the number of all possible partitions of $k$.

Examples of such a CW complex $M$ include the infinite product
$\prod_{i=1}^{\infty} K(\Z,2i)$ (with the weak topology) of
 Eilenberg-MacLane spaces and $BU=\lim_{q\to \infty} BU_q$.
 Although the homotopy type of these topological spaces are
 different, their corresponding Betti numbers coincide.
\end{remark}

Now we will compute Betti numbers of $C_{4,2} (\P^n)$ ($n\geq 9$)
 in a different way. Since $C_{p,2}(\P^n)-\sp^2(G(p+1,n+1))=Q_{p,n}$, we have
 $H_i(C_{p,2}(\P^n),\sp^2(G(p+1,n+1))\cong H^{BM}_i(Q_{p,n})$ for all $i$, where $H^{BM}_i$ denotes
 the Borel-Moore homology.
Let $A_{p,n}$ be the fiber bundle over  $G(p+2,n+1)$ whose fiber is the
space of all quadric hypersurfaces in $\P^{p+1}$ and
let $B_{p,n}$ be the fiber bundle over $G(p+2,n+1)$ whose fiber is the space of pairs of hyperplanes in $\P^{p+1}$. 
From the definition of $Q_{p,n}$, we have $H_i(A_{p,n},B_{p,n})\cong H^{BM}_i(Q_{p,n})$ for all $i$. In particular, 
\begin{equation}\label{eqn3}
H_i(C_{4,2}(\P^n),\sp^2(G(5,n+1))\cong H_i(A_{4,n},B_{4,n})
\end{equation} 
for $i\geq 0$ and $n\geq 9$.

\begin{lemma} \label{lemma2.5}
Let $A_{4,n},B_{4,n}$ be defined as above.
$$
\beta_i(\sp^2(G(5,n+1))=\left\{
\begin{array}{lll}
 1,1,3,5,11;&\hbox{for $i=0,2,4,6,8$}\\
 0; & \hbox{for $i$ odd.}
\end{array}
 \right.
$$

$$
\beta_i(A_{4,n})=\left\{
\begin{array}{lll}
 1,2,4,7,12; & \hbox{for $i=0,2,4,6,8$.}\\
 0; & \hbox{for $i$ odd.}
\end{array}
 \right.
$$

$$
\beta_i(B_{4,n})=\left\{
\begin{array}{lll}
 1,2,5,9,17; & \hbox{for $i=0,2,4,6,8$.}\\
 0; & \hbox{for $i$ odd.}
\end{array}
 \right.
$$
\end{lemma}

\bp To show the first formula, we note that all the odd Betti numbers of $G(5,n+1)$
are zero and the first five even Betti numbers of $G(5,n+1)$ are given by
$$\beta_i(G(5,n+1))=1,1,2,3,5~  \hbox{for $i=0,2,4,6,8$.}
$$
Therefore all the odd Betti numbers of $\sp^2(G(5,n+1))$
vanish and the first five even Betti numbers of $\sp^2(G(5,n+1))$ are given by (a special case of
MacDonald formula \cite{Macdonald})
$$
\beta_i(\sp^2(G(5,n+1)))=1,1,3,5,11 ~ \hbox{for $i=0,2,4,6,8$.}
$$

To show the second formula, we note that  $A_{4,n}$ is a fiber bundle over
$G(6,n+1)$ with fibers
the space of all quadric hypersurfaces in $\P^{5}$, i.e., fibers are isomorphic to $\P^{20}$.
By Lemma \ref{lemma2.1}, all the odd Betti numbers of $A_{4,n}$ vanish since both
$G(6,n+1)$ and $\P^{20}$
only carry non-vanishing even Betti numbers. Again, by Lemma \ref{lemma2.1},
\begin{equation}\label{eq8}
\beta_{2k}(A_{4,n})=\bigoplus_{i+j=k} \beta_{2i}(G(6,n+1))\cdot\beta_{2j}(\P^{20}).
\end{equation}
The first five even Betti numbers of $G(6,n+1)$ are given by
$$\beta_{i}(G(6,n+1))=1,1,2,3,5~ \hbox{for $i=0,2,4,6,8$.}$$
Hence from Equation (\ref{eq8}), we get the first five even Betti numbers of $\widetilde{Y}$:
\begin{equation}\label{eq9}
\beta_{i}(A_{4,n})=1,2,4,7,12 ~ \hbox{for $i=0,2,4,6,8$.}
\end{equation}

To show the third formula, we note  that $B_{4,n}$ is a fiber bundle over $G(6,n+1)$ with fibers the space of
pairs of hyperplanes in $\P^{p+1}$, i.e., fibers are isomorphic to $\sp^2 (\P^5)$.
By Lemma \ref{lemma2.1}, all the odd Betti numbers of $B_{4,n}$ vanish and the even Betti numbers of $B_{4,n}$ are given by the formula:
 \begin{equation}\label{eq10}
\beta_{2k}(B_{4,n})=\bigoplus_{i+j=k} \beta_{2i}(G(6,n+1))\cdot\beta_{2j}(\sp^2 (\P^5)).
\end{equation}
The first five Betti numbers of $\sp^2 (\P^5)$ are given as follows (cf. \cite{Macdonald}):
$$
\beta_{i}(\sp^2 (\P^5))=1,1,2,2,3~ \hbox{for $i=0,2,4,6,8$.}
$$
Therefore the five Betti numbers of $Z$ are given by the formula:
 \begin{equation}
\beta_{i}(B_{4,n})=1,2,5,9,17~ \hbox{for $i=0,2,4,6,8$.}
\end{equation}
\qe

\begin{proposition}\label{prop2.3}
The relations among the first 10 Betti numbers of
$C_{4,2}(\P^n)$ ($n\geq 9$) are given as follows:
$$
\beta_{2i}(C_{4,2}(\P^n))-\beta_{2i+1}(C_{p,2}(\P^n))=
1,1,2,3,6 \quad\hbox{for $i=0,1,2,3,4.$}
$$
In particular, $\beta_{8}(C_{4,2}(\P^n))\geq 6$.
\end{proposition}

\bp Set $M=C_{4,2}(\P^n)$ and $X=\sp^2 G(5,n+1)$. From the long exact sequence of homology groups for the pair $(M,X)$, we have
\begin{equation}\label{eqn8}
\cdots \to H_i(X)\to H_i(M)\to H_i(M,X)\to H_{i-1}(X)\to\cdots.
\end{equation}
Since $H_{2i-1}(X)=0$ for all $i$, Equation (\ref{eqn8}) breaks into exact sequences
\begin{equation}\label{eqn9}
0 \to H_{2i+1}(M)\to H_{2i+1}(M,X)\to H_{2i}(X)\to H_{2i}(M)\to H_{2i}(M,X)\to 0.
\end{equation}

Set $Y=A_{4,n}$ and $Z=B_{4,n}$.  From the long exact sequence of homology groups for the pair $(Y,Z)$, we have
\begin{equation}\label{eqn10}
 \cdots \to H_i(Z)\to H_i(Y)\to H_i(Y,Z)\to H_{i-1}(Z)\to\cdots
\end{equation}
Since $H_{2i-1}(Y)=0$ and $H_{2i-1}(Z)=0$ for all $i$, Equation (\ref{eqn10}) breaks into  exact sequences
\begin{equation}\label{eqn11}
0 \to H_{2i+1}(Y,Z)\to H_{2i}(Z)\to H_{2i}(Y)\to H_{2i}(Y,Z)\to 0.
\end{equation}

From equation (\ref{eqn3}),(\ref{eqn9}) and (\ref{eqn11}), we have
\begin{equation*}
\beta_{2i+1}(M)-\beta_{2i}(Z)+\beta_{2i}(Y)+\beta_{2i}(X)-\beta_{2i}(M)=0
\end{equation*}
i.e.,
\begin{equation}\label{eqn12}
\beta_{2i+1}(C_{4,2}(\P^n))-\beta_{2i}(B_{4,n})+\beta_{2i}(A_{4,n})+\beta_{2i}(\sp^2 G(5,n+1))-\beta_{2i}(C_{4,2}(\P^n))=0
\end{equation}

Now the proposition follows from equation (\ref{eqn12}) and Lemma \ref{lemma2.5}.

\qe


The contradiction between Corollary \ref{cor2.1}  and Proposition \ref{prop2.3} comes from the assumption
that the answer to Question \ref{ques1} for $d=2$ is affirmative. Therefore the answer to Question \ref{ques1} for
$d=2$ is negative, i.e.,  the induced map
$\pi_*(\D(d))\to \pi_*(\D(\infty))$  by inclusion is not always injective for $d=2$.
This completes the proof of Theorem \ref{Th1.2}.

\begin{remark}
We actually used the assumption
that $\pi_*(\D(2))\cong \pi_*(\D(\infty))$ are isomorphisms as abstract groups for $k\leq 9$ in the proof of Theorem \ref{Th1.2}. Hence
$\pi_*(\D(2))$ is not  isomorphic  to $\pi_*(\D(\infty))$ for all $*$ as abstract abelian groups.
\end{remark}


\section{Proof of Theorem \ref{prop1.1}}\label{sec:3}
In this section we will prove Theorem \ref{prop1.1}. The method comes from Lawson's  proof to
the Complex Suspension Theorem \cite{Lawson1}, i.e., the complex suspension to the space of $p$-cycles yields a homotopy equivalence to the
space of $(p+1)$-cycles. Here we briefly review the general construction of such a homotopy equivalence. For details, the reader  is referred to
\cite{Lawson1}, \cite{Lawson2} and \cite{Friedlander1}.

Fix a hyperplane $\P^n\subset\P^{n+1}$ and a point $\P^0\in \P^{n+1}-\P^{n}$.
For any non-negative integer $p$ and $d$, set
 $$T_{p+1,d}(\P^{n+1}):=\big\{c=\sum n_iV_i\in C_{p+1,d}(\P^{n+1})|\dim(V_i\cap\P^n)=p, ~\forall i\big \}.$$
(when $d=0$, $C_{p,0}(\P^n)$ is defined to be the empty cycle.)
 \begin{proposition}[\cite{Lawson1}]
The set  $T_{p+1,d}(\P^{n+1})$  is Zariski open in $C_{p+1,d}(\P^{n+1})$. Moreover, $T_{p+1,d}(\P^{n+1})$ is homotopy equivalent to
$C_{p,d}(\P^{n})$. In particular, their corresponding homotopy groups are isomorphic, i.e.,
\begin{equation}\label{eq12}
\pi_*(T_{p+1,d}(\P^{n+1}))\cong\pi_*(C_{p,d}(\P^{n})).
\end{equation}
\end{proposition}

\medskip
Fix linear embedding $\P^{n+1}\subset \P^{n+2}$ and two points $x_0,x_1\in \P^{n+2}-\P^{n+1}$.
Each projection $p_i:\P^{n+2}-\{x_0\}\to \P^{n+1}$ ($i=0,1$) gives us a holomorphic line bundle over $\P^{n+1}$.

Let $D\in C_{n+1,e}(\P^{n+2})$ be effective divisor of degree $e$ in $\P^{n+2}$
such that $x_0, x_1$ are not in $D$. Any effective cycle $c\in C_{p+1,d}(\P^{n+1})$
 can be lifted to a cycle with support in $D$, defined as follows:
$$ \Psi_D(c)=(\Sigma_{x_0}c)\cdot D.
$$

 The map $\Psi(c, D):=\Psi_D$ is a continuous map with variables $c$ and $D$.
Hence we have a continuous map $\Phi_D:C_{p+1,d}(\P^{n+1})\to C_{p+1,de}(\P^{n+2}-\{x_0,x_1\})$.
The composition of $\Phi_D$ with the projection $(p_0)_*$ is
$(p_0)_*\circ\Phi_D=e$ (multiplication by the integer $e$ in the monoid, $e\cdot c=c+\cdots+c$ for $e$ times).
The composition of $\Phi_D$ with the projection $(p_1)_*$ gives us a transformation of cycles in
 $\P^{n+1}$ which makes most of them intersecting properly to $\P^{n}$.
To see this, we consider the family of divisors $tD$, $0\leq t\leq 1$, given by scalar multiplication by $t$ in the line bundle
$p_0:\P^{n+2}-\{x_0\}\to \P^{n+1}$.

Assume $x_1$ is not in $tD$ for all $t$. Then the above construction gives us a family transformation
$$ F_{tD}:=(p_1)_*\circ \Psi_{tD}: C_{p+1,d}(\P^{n+1})\to C_{p+1,de}(\P^{n+1})
$$
for $0\leq t\leq 1$. Note that $F_{0D}\equiv d$(multiplication by $d$).

The question is that for a fixed $c$, which divisors $D\in C_{n+1,e}(\P^{n+2})$
($x_0$ is not in $D$ and $x_1$ is not in $\bigcup_{0\leq t\leq 1} tD$) have the property that
$$ F_{tD}(c)\in T_{p+1,de}(\P^{n+1})
$$
for all $0< t\leq 1$.

Set $B_c:=\{D\in C_{n+1,e}(\P^{n+2})|F_{tD}(c)~\hbox{is not in} ~T_{p+1,de} (\P^{n+1})~ \hbox{ for some }~ 0<t\leq 1\}$,
i.e., all degree $e$ divisors on $\P^{n+2}$
 such that some component of
 \begin{equation*}
 (p_1)_*\circ \Psi_{tD}(c)\subset \P^n
 \end{equation*}
 for some $t>0$.

\begin{proposition}[\cite{Lawson1}]
For $c\in C_{p+1,d}(\P^{n+1})$, ${\rm codim}_{\C}B_c\geq \big(^{p+e+1}_{~~~e}\big)$.
\end{proposition}

In this construction, if we take $e=1$, then $F_{tD}$ maps $C_{p+1,d}(\P^{n+1})$ to itself, i.e.,
$$ F_{tD}:=(p_1)_*\circ \Psi_{tD}: C_{p+1,d}(\P^{n+1})\to C_{p+1,d}(\P^{n+1}).
$$
Moreover, the image of $F_{tD}$ is in the Zariski open subset $T_{p+1,d}(\P^{n+1})$ if $D$ is not $B_c$.
We can find such a $D$ if ${\rm codim}_{\C}B_c\geq \big(^{p+1+1}_{~~~1}\big)=p+2$ is positive.

Suppose now that $f:S^k\to C_{p+1,d}(\P^{n+1})$ is a continuous map for $0<k\leq 2p+2$.
We may assume that $f$ is piecewise linear up to homotopy. Then the map $f$ is homotopic
to a map $S^k\to T_{p+1,d}(\P^{n+1})$. To see this, we consider the family
$$F_{tD}\circ f:S^k\to C_{p+1,d}(\P^{n+1}), \quad 0\leq t\leq 1,
$$
where $D$ lies outside the union $\bigcup_{x\in S^k}B_{f(x)}$. This is a set of real
codimension bigger than or equal to $2(p+2)-(k+1)$. Therefore, $2(p+2)-(k+1)\geq 1$, i.e.,
 $k\leq 2p+2$, then such a $D$ exists. This proves that the map $i_*:\pi_k(T_{p+1,d}(\P^{n+1}))\to \pi_k(C_{p+1,d}(\P^{n+1}))$
induced by inclusion $i:T_{p+1,d}(\P^{n+1}))\hookrightarrow C_{p+1,d}(\P^{n+1})$
is \emph{surjective} if $k\leq 2p+2$.

Similarly, suppose that $g:(D^{k+1},S^k)\to (C_{p+1,d}(\P^{n+1}),T_{p+1,d}(\P^{n+1}))$
is a pair of continuous map for $0<k\leq 2p+1$. Then the map $g$ can be deformed through a map of pairs to
$\tilde{g}:(D^{k+1},S^k)\to (T_{p+1,d}(\P^{n+1}),T_{p+1,d}(\P^{n+1}))$
if $2(p+2)-(k+2)\geq 1$, i.e., $k\leq 2p+1$. This proves that the map $i_*:\pi_k(T_{p+1,d}(\P^{n+1}))\to \pi_k(C_{p+1,d}(\P^{n+1}))$
induced by inclusion $i:T_{p+1,d}(\P^{n+1}))\hookrightarrow C_{p+1,d}(\P^{n+1})$
is \emph{injective} if $k\leq 2p+1$.

Therefore,
\begin{equation}\label{eq13}
\pi_k(T_{p+1,d}(\P^{n+1}))\stackrel{i_*}{\cong} \pi_k(C_{p+1,d}(\P^{n+1}))
\end{equation}
for $0\leq k\leq 2p+1$.

The combination of Equation (\ref{eq12}) and (\ref{eq13}) gives us the following result:
\begin{proposition}\label{prop3.3}
The complex suspension $\Sigma:C_{p,d}(\P^{n})\to C_{p+1,d}(\P^{n+1})$ induces an isomorphism
\begin{equation}\label{eq14}
 \Sigma_*:\pi_k(C_{p,d}(\P^{n})){\cong} \pi_k(C_{p+1,d}(\P^{n+1}))
\end{equation}
 for $0\leq k\leq 2p+1$.
\end{proposition}

As a corollary,  we get the simply connectedness of $C_{p,d}(\P^n)$, which has been obtained using general position arguments
by Lawson (\cite{Lawson1}, the proof to Lemma 2.6.):
\begin{corollary}[\cite{Lawson1}]\label{cor3.1}
The Chow variety $C_{p,d}(\P^n)$ is simply connected for  integers $p,d,n\geq 0$.
\end{corollary}
\bp Since $C_{0,d}(\P^{n})$ can be identified with $d$-th symmetric product $\sp^d(\P^n)$ of $\P^n$ and $\sp^d(\P^n)$
is path connected, we have $\pi_0(C_{0,d}(\P^{n}))=0$ for all $d,n\geq 0$. Repeating using Equation (\ref{eq14}), we know
$\pi_0(C_{p,d}(\P^{n}))=0$ for all $p,d,n\geq 0$.
Moreover, since $\sp^d(\P^n)$ is simply connected for all $d,n\geq 0$, we have $\pi_1(C_{0,d}(\P^{n}))=0$ for all $d,n\geq 0$.
Repeating using Equation (\ref{eq14}), we get $$\pi_1(C_{p,d}(\P^{n}))\cong \pi_1(C_{p-1,d}(\P^{n-1}))\cong\cdots\cong\pi_1(C_{0,d}(\P^{n-p}))=0$$
 for all $p,d,n\geq 0$.

\qe

Now we study the connectedness of maps induced by the inclusion $i:\P^n\hookrightarrow \P^{n+1}$.

\begin{proposition}\label{prop3.5}
 For any integer $d\geq 1$, the inclusion $i:C_{p,d}(\P^{n})\hookrightarrow C_{p,d}(\P^{n+1})$ induces an isomorphism
\begin{equation}\label{eq15}
 \pi_k(C_{p,d}(\P^{n}))\stackrel{i_*}{\cong} \pi_k(C_{p,d}(\P^{n+1}))
\end{equation}
for $0\leq k\leq 2(n-p)$.
\end{proposition}

\begin{remark}\label{remark3.6}
By using Proposition \ref{prop3.5},  we give another possibly more elementary proof Corollary \ref{cor3.1}. If
$n=p$, then $C_{p,d}(\P^{n})$ is a point and so it is simply connected. If $n=p+1$, then $C_{p,d}(\P^{n})\cong \P^{(^{n+d}_{~d})-1}$
so it is simply connected. If $n-p\geq 2$, then $\pi_k(C_{p,d}(\P^{n})){\cong} \pi_k(C_{p,d}(\P^{n-1}))\cong\cdots\cong \pi_k(C_{p,d}(\P^{p+1}))=0$
for $k\leq 1$ by using Proposition \ref{prop3.5} and so $C_{p,d}(\P^{n})$ is simply connected.
\end{remark}

Proposition \ref{prop3.5} can be used to compute the second homotopy group of Chow varieties.

\begin{corollary}
For $d\geq 1$ and $n>p \geq 0$, we have $\pi_2(C_{p,d}(\P^{n})){\cong}\Z$ and  $H_2(C_{p,d}(\P^{n})){\cong}\Z$.
\end{corollary}
\bp  Replacing $\pi_k$ by $\pi_2$ in Remark \ref{remark3.6} yields the proof of the first statement. The second statement is a result of the first statement, Corollary \ref{cor3.1} and the Hurewicz isomorphism theorem.

\qe

Lawson's idea in the proof of the Complex Suspension Theorem in \cite{Lawson1} can be used to prove
Proposition \ref{prop3.5}.

For any non-negative integer $p$ and $d$, set
$$U_{p,d}(\P^{n+1}):=\big\{c=\sum n_iV_i\in C_{p,d}(\P^{n+1})|\hbox{$\P^0$ is not in $\cup_i V_i$}\big\}.$$

Proposition \ref{prop3.5} follows directly from the  combination of  Lemma \ref{lemma3.1} and \ref{lemma3.2} below:

\begin{lemma}\label{lemma3.1}
$U_{p,d}(\P^{n+1})$ is homotopy equivalent to
$C_{p,d}(\P^{n})$. In particular, their corresponding homotopy groups are isomorphic, i.e.,
\begin{equation}\label{eq16}
\pi_*(U_{p,d}(\P^{n+1}))\cong\pi_*(C_{p,d}(\P^{n})).
\end{equation}
\end{lemma}

\bp Let $p_0:\P^{n+1}-\P^0\to \P^n$ be the canonical projection away from $\P^0\in \P^{n+1}-\P^n$. Then
$p_0$ induces a deformation retract from $U_{p,d}(\P^{n+1})$  to $C_{p,d}(\P^{n})$.

To see this, note that $p_0$ is a holomorphic line bundle and
let
$F_t:(\P^{n+1}-\P^0)\times \C\to \P^{n+1}-\P^0 $ denote the scalar multiplication by $t\in \C$ in this bundle.
This map $F_t$ is holomorphic (in fact, algebraic) and satisfies $F_1=id_{\P^{n+1}-\P^0}$
and $F_0=p_0$. Hence $F_t$ induces a family of continuous maps
$(F_t)_*: U_{p,d}(\P^{n+1}) \to C_{p,d}(\P^{n})$. Therefore, $(p_0)_*$ is a deformation retraction.

\qe

\begin{lemma}\label{lemma3.2}
The inclusion $i: U_{p,d}(\P^{n+1})\hookrightarrow C_{p,d}(\P^{n+1})$ is $2(n-p)$-connected.
\end{lemma}

\bp By definition, it is enough to show that the induced maps on homotopy groups
$$i_*: \pi_k(U_{p,d}(\P^{n+1}))\to \pi_k(C_{p,d}(\P^{n+1}))$$ are isomorphisms for
$k\leq 2(n-p)$. Let $f:S^k\to C_{p,d}(\P^{n+1})$ be a continuous map for $k\leq 2(n-p)$.
 We may assume $f$ to be piecewise linear up to homotopy. Then $f$ is
 homotopy to a map $S^k\to U_{p,d}(\P^{n+1})$. To see this, we note firstly that the union
$$
\bigcup_{x\in S^k} f(x)
$$
is a set of real codimension $\geq 2(n+1)-2p-k\geq 2>0$. So we can find a point $Q\in \P^{n+1}-\P^{n}$ such
that $Q$ is not in $\bigcup_{x\in S^k} f(x)$. Let $G_t$ be a family of automorphism
of $\P^{n+1}$ mapping $\P^0$ to $Q$ but preserving $\P^n$.
Composing with the automorphism $G_t$, we obtain the family $G_t\circ f: S^k \to C_{p,d}(\P^{n+1})$
such that $G_0\circ f=f$ and $G_1\circ f:S^k\to U_{p,d}(\P^{n+1})$.
Hence $i_*$ is \emph{surjective} for $k\leq 2(n-p)$.

Similarly, suppose $g$ is a map of pairs $g:(D^{k+1}, S^k)\to (C_{p,d}(\P^{n+1}), U_{p,d}(\P^{n+1}))$.
Then the map can be deformed through a map of pairs to one with image in $U_{p,d}(\P^{n+1})$ if $k\leq 2(n-p)$.
Therefore, $i_*$ is \emph{injective} for $k\leq 2(n-p)$.

\qe

By Proposition \ref{prop3.5}, $\pi_k(C_{p,d}(\P^{n}))$ is stable when $n\to \infty$.
By the combination of Equations (\ref{eq14}) and (\ref{eq15}), we have  the following isomorphism
$$\pi_k(C_{p,d}(\P^{n}))\cong \lim_{m,q\to\infty}\pi_k(C_{p+q,d}(\P^{n+m+q})
$$
for $0\leq k\leq 2p+1$ and $k\leq 2(n-p)$. This completes the proof of Theorem \ref{prop1.1}.
\qe

Department of Mathematics, Massachusetts Institute of Technology, Room 2-363B,
 77 Massachusetts Avenue,
 Cambridge, MA 02139, USA

Email: {wenchuan@math.mit.edu}

\end{document}